\documentclass[12pt]{article}

\usepackage{amsmath,amsthm,amsfonts,amssymb,amscd,cite}
\usepackage{tikz}
\usepackage{color}
\usepackage{hyperref}

\newtheorem{claim}{Claim}
\newtheorem{theorem}{Theorem}[section]

\newtheorem{corollary}[theorem]{Corollary}

\newtheorem{problem}[theorem]{Problem}

\date{}
\begin{document}

\title{On the Anti-Ramsey Number Under Edge Deletion}

\author{
Ali Ghalavand$^{a,}$\thanks{ email:\texttt{alighalavand@nankai.edu.cn }} \and  Qing Jie$^{b,}$\thanks{ email:\texttt{qingjie1126@163.com}} \and  Zemin Jin$^{b,}$\thanks{Corresponding author, email:\texttt{zeminjin@zjnu.cn }}\and  Xueliang Li$^{a,c,}$\thanks{ email:\texttt{lxl@nankai.edu.cn }}\and Linshu  Pan $^{b,}$\thanks{ email:\texttt{1925984948@qq.com }} }

\maketitle

\begin{center}
$^a$ Center for Combinatorics and LPMC, Nankai University, Tianjin 300071, China \\
$^b$ School of Mathematical Sciences, Zhejiang Normal University
Jinhua 321004, China\\
$^c$ School of Mathematical Sciences, Xinjiang Normal University, Urumchi, Xinjiang, China
\end{center}

\begin{abstract}
 According to a study by Erdős et al. in 1975, the anti-Ramsey number of a graph \(G\), denoted as \(AR(n, G)\), is defined as the maximum number of colors that can be used in an edge-coloring of the complete graph \(K_n\) without creating a rainbow copy of \(G\). 
In this paper, we investigate the anti-Ramsey number under edge deletion and demonstrate that both decreasing and unchanging are possible outcomes. For three non-negative integers \(k\), \(t\), and \(n\), let \(G = kP_4 \cup tP_2\). Let \(E'\) be a subset of the edge set \(E(G)\) such that every endpoint of these edges has a degree of two in \(G\). 
We prove that if one of the conditions
(i) \(t \geq k + 1 \geq 2\) and \(n \geq 8k + 2t - 4\); 
(ii) \(k, t \geq 1\) and \(n = 4k + 2t\); 
(iii) \(k = 1\), \(t \geq 1\), and \(n \geq 2t + 4\), occurs then the behavior of the anti-Ramsey number remains consistent when the edges in \(E'\) are removed from \(G\), i.e., \(AR(n, G) = AR(n, G - E')\). However, this is not the case when \(k \geq 2\), \(t = 0\), and \(n=4k\).
As a result, we calculate \(AR(n,kP_4 \cup tP_2)\) for the cases:
(i) \(t \geq k + 1 \geq 2\) and \(n \geq 8k + 2t - 4\); 
(ii) \(k, t \geq 1\) and \(n = 4k + 2t\); 
(iii) \(k = 1\), \(t \geq 0\), and \(n \geq 2t + 4\); 
(iv) \(k \geq 1\), \(t = 0\), and \(n = 4k\).     
\end{abstract}


\section{Introduction}
This study examines the impact of edge deletion on anti-Ramsey numbers for linear forests made up of paths of lengths two and four (i.e., \( P_2 \) and \( P_4 \)). We will start by reviewing important definitions and providing some historical context.
Let \( G \) be a finite, simple graph. An {\it edge-coloring} of \( G \) is a function \( c : E(G) \to \mathcal{C} \), where \( \mathcal{C} \) is a set of colors. A subgraph \( H \subseteq G \) is considered {\it rainbow} if all its edges have distinct colors according to \( c \). For a subgraph \( H \) or an edge set \( E' \subseteq E(G) \), we denote the set of colors on its edges by \( c(H) \) or \( c(E') \), respectively.
For a vertex subset \( V_1 \subseteq V(G) \), the {\it induced subgraph} \( G[V_1] \) has the vertex set \( V_1 \) and includes every edge of \( G \) whose endpoints are in \( V_1 \). For another subset \( V_2 \) of \( V(G) \), the notation \( [V_1, V_2]_G \) represents the set of edges in \( G \) that have one endpoint in \( V_1 \) and the other in \( V_2 \). The {\it disjoint union} of two graphs \( G \) and \( G' \) is denoted \( G \cup G' \). 
For a positive integer \( k \), \( kG \) denotes the disjoint union of \( k \) copies of \( G \). For a subset \( E' \) of \( E(G) \), the notation \( G - E \) denotes the graph obtained from \( G \) by removing the edges in \( E' \).
 For a positive integer \( z \), we define \( [z] = \{1, 2, \ldots, z\} \) and \( [0] = \emptyset \). We denote the complete graph, cycle, and path on \( l \) vertices by \( K_l \), \( C_l \), and \( P_l \), respectively.

For a positive integer \( n \) and a graph \( H \), the {\it anti-Ramsey number} \( AR(n, H) \) is defined as the maximum number of colors in an edge-coloring of the complete graph \( K_n \) that does not contain a rainbow copy of \( H \). This concept was introduced by Erdős et al. in 1975 \cite{Er1}, where they established its connection to Turán numbers and conjectured values for the anti-Ramsey numbers of paths and cycles.

Subsequent research has extensively investigated anti-Ramsey numbers for edge-disjoint subgraphs, as well as for graphs composed of multiple small, disjoint components. Results concerning the anti-Ramsey numbers for spanning trees in complete graphs and more general settings can be found in \cite{R2-SP-1,R2-SP-2,p1-9,R2-SP-3,R2-SP-4}. The anti-Ramsey numbers for matchings (collections of disjoint edges) in complete graphs are discussed in \cite{p1-2, p1-5, p1-7, p1-9, p1-14}, while extensions to uniform hypergraphs appear in \cite{R2-HP-1,R2-HP-2,R2-HP-3}.
More recently, research has shifted toward more intricate structures, particularly vertex-disjoint unions of small connected graphs and matchings. Our work advances this line of inquiry by examining linear forests composed of short odd paths. For clarity, we summarize the existing results most relevant to our study below. Simonovits and Sós \cite{p2-10} proved that if \( t \geq 2 \) and \( n \) is sufficiently large, then the following holds:
\[
AR(n,P_t) = {\left\lfloor\frac{t-1}{2}\right\rfloor\choose2} + \left( \left\lfloor\frac{t-1}{2}\right\rfloor - 1 \right) \left(n - \left\lfloor\frac{t-1}{2}\right\rfloor + 1\right) + 1 + \epsilon,
\]
where \( \epsilon = 1 \) if \( t \) is even and \( \epsilon = 0 \) otherwise. 
Montellano-Ballesteros and Neumann-Lara \cite{p3-6} established that for any \( n \geq t \geq 3 \),
\[
AR(n,C_t) = {t-1\choose2} \left\lfloor\frac{n}{t-1}\right\rfloor + \left\lceil\frac{n}{t-1}\right\rceil + {n \mod (t-1) \choose 2}.
\]
Schiermeyer \cite{p2-9} demonstrated that for all \( t \geq 2 \) and \( n \geq 3t + 3 \), the following holds:
\[
AR(n,tP_2) = {t-2 \choose 2} + (t-2)(n - t + 2) + 1.
\]
Chen et al. \cite{p1-2} and Fujita et al. \cite{p1-5} showed that for all \( t \geq 2 \) and \( n \geq 2t + 1 \):
\begin{equation}\label{eqth0}
AR(n,tP_2) = \begin{cases} 
(t-2)(2t-3) + 1 & \text{\rm when } n \leq \frac{5t-7}{2}, \\ 
(t-2)(n - \frac{t-1}{2}) + 1 & \text{ \rm when } n \geq \frac{5t-7}{2}. 
\end{cases}
\end{equation}
The remaining case \( n = 2t \) was addressed by Haas and Young \cite{p1-7}, who confirmed the conjecture made in \cite{p1-5}. They showed that if \( n = 2t \), then:
\begin{equation}\label{eqth1}
AR(n,tP_2) = \begin{cases} 
3&\text{\rm when } t=2,\\
\frac{1}{2}(t-2)(3t+1) + 1 & \text{\rm when } 3 \leq t \leq 6, \\ 
(t-2)(2t-3) + 2 & \text{\rm when } t \geq 7. 
\end{cases}
\end{equation}
Bialostocki et al. \cite{p1-1} determined the anti-Ramsey number for the graphs \( P_4 \) and \( P_4 \cup P_2 \). They established the following results:
\begin{equation}\label{eqth2}
AR(n, P_4 \cup tP_2) = \begin{cases} 
3 & \text{when } t = 0 \text{ and } n = 4, \\ 
2 & \text{when } t = 0 \text{ and } n \geq 5, \\ 
n & \text{when } t = 1 \text{ and } n \geq 6.
\end{cases}
\end{equation}
They also proved that \( AR(n,P_3 \cup P_2) = 2 \) when \( n \geq 5 \) and \( AR(n,P_3 \cup 2P_2) = n \) when \( n \geq 7 \).
For more general classes, Gilboa and Roditty \cite{p1-6} derived upper bounds on the anti-Ramsey numbers for graphs of the form \( L \cup tP_2 \) and \( L \cup tP_3 \) under specific conditions. They showed that for sufficiently large $n$:
\begin{enumerate}
 \item[] $AR(n,P_{k+1}\cup tP_3) = (t+\lfloor\frac{k}{2}\rfloor-1)(n-\frac{t+\lfloor\frac{k}{2}\rfloor}{2}) + 1 + (k \mod 2)$ when $k\geq3$;
 \item[] $AR(n,kP_3\cup tP_2) = (k+t-2)(n-\frac{t+k-1}{2}) + 1$ when $k,t\geq2$;
 \item[] $AR(n,P_2\cup tP_3) = (t-1)(n-\frac{t}{2}) + 1$ when $t\geq1$;
 \item[] $AR(n,P_3\cup tP_2) = (t-1)(n-\frac{t}{2}) + 1$ when $t\geq2$;
 \item[] $AR(n,P_4\cup tP_2) = t(n-\frac{t+1}{2}) + 1$ when $t\geq1$;
 \item[] $AR(n,C_3\cup tP_2) = t(n-\frac{t+1}{2}) + 1$ when $t\geq1$;
 \item[] $AR(n,tP_3) = (t-1)(n-\frac{t}{2}) + 1$ when $t\geq1$.
\end{enumerate}
He and Jin \cite{p1-8} established that if \( t \geq 2 \) and \( n \geq 2t + 3 \), then 
\[
AR(n, P_3 \cup tP_2) = 
\begin{cases} 
t(2t-1) + 1 & \text{when } 2t + 3 \leq n \leq \left\lfloor \frac{5t + 2}{2} + \frac{1}{t-1} \right\rfloor, \\
\frac{1}{2} (t-1)(2n - t) + 1 & \text{when } n \geq \left\lceil \frac{5t + 2}{2} + \frac{1}{t-1} \right\rceil. 
\end{cases}
\]
They also proved that if \( t \geq 2 \) and \( n \geq 2t + 7 \), then 
\[
AR(n, 2P_3 \cup tP_2) = 
\begin{cases} 
(t+1)(2t + 3) + 1 & \text{when } 2t + 7 \leq n \leq \left\lfloor \frac{5t + 11}{2} + \frac{3}{t} \right\rfloor, \\
\frac{1}{2} t(2n - t - 1) + 1 & \text{when } n \geq \left\lceil \frac{5t + 11}{2} + \frac{3}{t} \right\rceil. 
\end{cases}
\]
Let \( \lambda = \frac{9k + 5t - 7}{2} + \frac{k(k+1)}{2(k+t-2)} \). Jie et al. \cite{p1-10} studied the anti-Ramsey number of \( kP_3 \cup tP_2 \) and proved that if \( k \geq 2 \), \( t \geq \frac{k^2 - k + 4}{2} \), and \( n \geq 3k + 2t + 1 \), then 
\[
AR(n, kP_3 \cup tP_2) = 
\begin{cases} 
\frac{1}{2} (3k + 2t - 3)(3k + 2t - 4) + 1 & {\rm~when~ } n \leq \left\lfloor \lambda \right\rfloor, \\
\frac{1}{2} (k + t - 2)(2n - k - t + 1) + 1 & {\rm~when~ } n \geq \left\lceil \lambda \right\rceil. 
\end{cases}
\] 
 Recently, two of the present authors \cite{new-2} investigated the anti-Ramsey numbers for spanning linear forests composed of paths of lengths 2 and 3. They demonstrated that for any integers \( k \), \( t \), and \( n \) satisfying the conditions \( k \geq 2 \), \( t \geq \frac{k^2 - 3k + 4}{2} \), and \( n = 2t + 3k \), the following holds:
\[
AR(n, kP_3 \cup tP_2) = (3k + 2t - 3)\left(t + \frac{3k}{2} - 2\right) + 1.
\]
Building on this work, two other authors from the current team \cite{GL-1} showed that this result is applicable in more general cases. They proved that for any integers \( k \), \( t \), and \( n \) where \( k \geq 1 \), \( t \geq 2 \), and \( n = 2t + 3k \), the same equation holds:
\[
AR(n, kP_3 \cup tP_2) = (3k + 2t - 3)\left(t + \frac{3k}{2} - 2\right) + 1.
\]
Jin and Gu \cite{new-1} investigated the anti-Ramsey number of graphs where each component is either \( K_4 \) or \( P_2 \). They established that for any integers \( n \geq \max\{7, 2t + 4\} \) and \( t \geq 1 \), the following holds:
\[
AR(n, K_4 \cup tP_2) =
\begin{cases}
\lfloor \frac{n}{2} \rfloor \lceil \frac{n}{2} \rceil + 1, & \text{when } t - 1 \leq t_n, \\
t(n - t) + \binom{t}{2} + 1, & \text{when } t - 1 > t_n.
\end{cases}
\]
Here, \( t_n \) is a constant derived from a specific function of \( n \).
Additionally, the authors in \cite{new-3} studied the anti-Ramsey number of linear forests consisting of paths of lengths 2 and 5. They proved that if \( n_0(t) = \frac{5t + 7}{2} + \frac{1}{t} \), then the result is as follows:
\[
AR(n, P_5 \cup tP_2) =
\begin{cases}
n + 1, & \text{when } t = 1 \text{ and } n \geq 7, \\
t(n - t) + \binom{t}{2} + 1, & \text{when } t \leq 4 \text{ and } n \geq 2t + 6,\\
(t + 1)(2t + 1) + 1, & \text{when } t\ge 5 \text{ and } 2t + 6 \leq n \leq \lfloor n_0(t) \rfloor, \\
t(n - t) + \binom{t}{2} + 1, & \text{when } t\ge 5 \text{ and } n \geq \lceil n_0(t) \rceil.
\end{cases}
\] 
Fang et al. \cite{p1-4} considered \( F \) to be a linear forest with components of order \( p_1, p_2, \ldots, p_t \), where \( t \geq 2 \) and \( p_i \geq 2 \) for \( 1 \leq i \leq t \), with at least one \( p_i \) being even. They proved that for sufficiently large \( n \), 
\[
AR(n, F) = \left( \sum_{i \in [t]} n \left\lfloor \frac{p_i}{2} \right\rfloor - \epsilon \right) n + O(1),
\] 
where \( \epsilon = 1 \) if all \( p_i \) are odd, and \( \epsilon = 2 \) otherwise.
Following this, Xie et al. \cite{p1-14} established an exact expression for the anti-Ramsey numbers of linear forests containing even components. They proved that if at least one \( p_i \) is even, where \( t \geq 2 \) and \( p_i \geq 2 \) for all \( 1 \leq i \leq t \), then for sufficiently large \( n \), 
\[
AR(n, F) = \binom{\sum_{i \in [t]} \left\lfloor \frac{p_i}{2} \right\rfloor - 2}{2} + \left( \sum_{i \in [t]} \left\lfloor \frac{p_i}{2} \right\rfloor - 2 \right) \left( n - \sum_{i \in [t]} \left\lfloor \frac{p_i}{2} \right\rfloor + 2 \right) + 1 + \epsilon,
\]
where \( \epsilon = 1 \) if exactly one \( p_i \) is even, or \( \epsilon = 0 \) if at least two \( p_i \) are even. Note that this expression does not account for the case when \( n = \sum_{i \in [t]} p_i \). 

In this study, we will investigate the impact of edge deletion on anti-Ramsey numbers for linear forests composed of paths of lengths two and four. We will demonstrate that both decreasing and no changes in these numbers are possible. Let \( E' \) be a subset of \( E(G) \). It is evident that for any edge coloring of \( K_n \), the presence of a rainbow subgraph in the form of \( G \) implies the existence of a rainbow subgraph isomorphic to \( G - E' \). Therefore, we can conclude that:

\begin{equation}\label{eqth3}
AR(n, G) \geq AR(n, G-E').
\end{equation}

By applying Equations \eqref{eqth0}, \eqref{eqth1}, and \eqref{eqth2}, we can deduce the following relationships:
\[
AR(4, P_4) = AR(4, 2P_2) \quad \text{and} \quad AR(n, P_4) = AR(n, 2P_2) + 1, \quad \text{for } n \geq 5.
\]
These results indicate that in certain cases of edge deletion in a graph, it is possible to either decrease the anti-Ramsey numbers or maintain them at the same level. 

Next, we will focus on the anti-Ramsey number of the linear forest \( kP_4 \cup tP_2 \) and prove the following theorems:
\begin{theorem}\label{th1}
	For three integers $k$, \( t \), and \( n \), if  \( t\geq k+1\geq 2 \), and \( n \ge 8k+2t-4 \), then
	\[
	AR(n, kP_4 \cup tP_2) = AR(n, (2k+t)P_2).
	\]
\end{theorem}
 
\begin{theorem}\label{th2}
	For three integers $k$, \( t \), and \( n \), if  \( k\geq 1 \), \(t\geq1\), and \( n =4k+2t \), then
	\[
	AR(n, kP_4 \cup tP_2) = AR(n, (2k+t)P_2).
	\]
\end{theorem}

\begin{theorem}\label{th3}
	For two integers $t$ and \( n \), if  \( t\geq 1 \), and \( n \ge 2t+4 \), then
	\[
	AR(n, P_4 \cup tP_2) = AR(n, (t+2)P_2).
	\]
\end{theorem}

Let \( k \), \( t \), and \( n \) be three positive integers. Define the graph \( G \) as \( G = kP_4 \cup tP_2 \). Let \( E' \) be a subset of the edge set \( E(G) \) such that every endpoint of these edges has a degree of two in \( G \). Based on Theorems \ref{th1}, \ref{th2}, and \ref{th3}, We can observe that \( AR(n, G) = AR(n, G - E') \) under the conditions: (i) if \( t \geq k + 1 \) and \( n \geq 8k + 2t - 4 \); (ii) if \( n = 4k + 2t \); or (iii) if \( k = 1 \) and \( n \geq 2t + 4 \). These conditions ensure that the behavior of the anti-Ramsey number remains consistent when the edges in \( E' \) are removed from \( G \).

\begin{theorem}\label{kp4}
	For two integers $k$ and $n$, if $k\ge 2$ and  $n= 4k$, then
	\[
	AR(n,kP_4)=(2k-1)(4k-3) + 1.  
	\]
\end{theorem}

 By applying Equation \eqref{eqth1} and Theorem \ref{kp4}, we can deduce the following:
\[
AR(4k,kP_4) - AR(4k,2kP_2) = 
\begin{cases} 
2k^2 - 5k + 4 & \text{when } 2 \leq k \leq 3, \\ 
4k - 4 & \text{when } k \geq 4. 
\end{cases}
\]
This leads us to investigate whether the condition \( t \geq 1 \) ensures that \( AR(n, G) = AR(n, G - E') \), where \( G = kP_4 \cup tP_2 \) and \( E' \) represents the specified edge set mentioned earlier. The following corollaries directly result from Theorems \ref{th1}, \ref{th2}, and \ref{th3}, as well as Equations \eqref{eqth0} and \eqref{eqth1}. It is worth noting that Corollary \ref{co1c} corrects some of the results found in \cite{n-p4-1}. 

\begin{corollary}
For three integers $k$, \( t \), and \( n \), if  \( t\geq k+1\geq 2 \), and \( n \ge 8k+2t-4 \), then
\[
	AR(n,kP_4\cup tP_2) =
\begin{cases} 
14&\text{\rm when } ~k=1,t=2,\\
(2k + t - 2)(4k + 2t - 3) + 1 & \text{\rm when } n \leq \frac{5(2k+t)-7}{2}, \\ 
(2k + t - 2)\left(n - \frac{2k+t-1}{2}\right) + 1 & \text{\rm when } n \geq \frac{5(2k+t)-7}{2}. 
\end{cases}
\]
\end{corollary}

\begin{corollary} 
	For three integers $k$, \( t \), and \( n \), if  \( k,t \geq 1 \), and \( n = 4k+2t \), then
	\[
	AR(n,kP_4\cup tP_2) = \begin{cases} 
		\frac{1}{2}(2k+t-2)(6k+3t+1) + 1 & {\rm~when~ } 2k+t\leq6,\\ 
		(2k+t-2)(4k+2t-3) + 2 & {\rm~when~ } 2k+t\geq7. 
	\end{cases}
	\]
\end{corollary}

\begin{corollary} \label{co1c}
	For two integers \( t \) and \( n \), if  \( t \geq 0 \), and \( n \geq 2t+4 \), then
	\[
	AR(n,P_4\cup tP_2) = \begin{cases} 
		3& {\rm when~ }~ t=0,\\
		\frac{1}{2}t(3t+7) + 1 & {\rm when~ }~ 1 \leq t \leq 4, n=2t+4, \\ 
		t(2t+1) + 2 & { \rm when~ }~ 5\leq t, n=2t+4, \\
		t(2t+1) + 1 & {\rm when~ }~ 2t+5\leq n \leq \frac{5t+3}{2}, \\ 
		t(n - \frac{t+1}{2}) + 1 & {\rm when~ }~\frac{5t+3}{2}\leq n. 
	\end{cases}
	\]
\end{corollary}

\section{Proof of Theorem~\ref{th1}}
For three integers \( k \), \( t \), and \( n \), let \( t \geq k + 1 \geq 2 \) and \( n \geq 8k + 2t - 4 \). By employing Equation \eqref{eqth3}, we find that 
\[
AR(n, kP_4 \cup tP_2) \geq AR(n, (2k + t)P_2).
\] 
To prove our theorem, we need to examine the following inequality: 
\begin{equation}\label{eqth4}
AR(n, kP_4 \cup tP_2) \leq AR(n, (2k + t)P_2).
\end{equation}
 To prove this, we need to demonstrate that if \( c \) is an arbitrary edge-coloring of \( K_n \) such that \( |c(K_n)| = AR(n, (2k+t)P_2) + 1 \), then \( K_n \) contains a rainbow subgraph isomorphic to \( kP_4 \cup tP_2 \). We will employ induction on \( k \).

Assume \( k = 1 \). Since \( |c(K_n)| = AR(n, (t+2)P_2) + 1 \), it follows that there exists a rainbow subgraph \( H_2 \) of \( K_n \) such that \( H_2 \cong (t+2)P_2 \). Let \( G \) be a spanning rainbow subgraph of \( K_n \) that includes \( H_2 \). For the next part of the proof, let's denote the edges of \( H_2 \) as \( E(H_2) = \{y^i_1 y^i_2 : i \in [t+2]\} \) and define \( V' = V(K_n) - V(H_2) \). 

Now, for the sake of contradiction, assume that \( K_n \) does not contain a rainbow \( P_4 \cup tP_2 \). We will present two claims based on this assumption.

\begin{claim}\label{cl1-p2-p2}
	$[V(P_2^i), V(P_2^j)]_G=\emptyset$ for $i\neq j\in [t+2]$.
\end{claim}
Suppose that there are two integers $i\neq j\in [t+2]$ such that $[V(P_2^i), V(P_2^j)]_G\neq \emptyset$. Without loss of generality, let $y_1^iy_1^j\in E(G)$. Then $y_2^iy_1^iy_1^jy_2^j$ and $H-V(P_2^i)-V(P_2^j)$ form a rainbow $P_4\cup tP_2$, a contradiction. Hence the claim holds. \qed

From Claim \ref{cl1-p2-p2}, we can derive 
\begin{equation}\label{eq-p2-p2}
	\begin{aligned} 
		|E(G[V(H_2))]| = t+2.
	\end{aligned}  
\end{equation}

\begin{claim}\label{cl2-2n-4t-8}
It holds that 
\begin{equation}\label{eq-2n-4t-8}
|E(G[V'])|+|[V(H_2),V']_G|\leq2n - 4t - 8.
\end{equation}
\end{claim}
Based on our assumptions, we find that any component of \( G[V'] \) that is not a star is isomorphic to either \( K_1 \), \( K_2 \), or \( K_3 \). Additionally, for indices \( i \neq j \in [t+2] \) and vertices \( u \neq v \in V' \), we can assert that the subgraph \( G[\{y^i_1, y^i_2, y^j_1, y^j_2, u\}] \) contains at least one isolated edge. Furthermore, for indices \( i \in [t+2] \), \( j \in [2] \), and a vertex \( v \in V' \), if \( y^i_j v \in E(G) \), then there does not exist a vertex \( u \in (V' - \{v\}) \) such that \( uv \in E(G) \).
As a result, we conclude that the number of edges in \( G \) with at least one endpoint in \( V' \) is at most \( 2|V'| = 2n - 4t - 8 \). Therefore, we have established the claim. \qed

By utilizing Equations \eqref{eq-p2-p2} and \eqref{eq-2n-4t-8}, we can conclude that 
\[
|c(K_n)| = |E(G)| \leq 2n - 3t - 6.
\]
However, applying Equations \eqref{eqth0} and \eqref{eqth1} leads us to the inequality \( |c(K_n)| = AR(n, (t + 2)P_2) + 1>2n - 3t - 6 \), which creates a contradiction.  Therefore, for \( k = 1 \), we affirm that Equation \eqref{eqth4} holds true.

Let \( k \geq 2 \). Given that 
\[
|c(K_n)| = AR(n, (2k + t) P_2) + 1 = AR(n, (2(k - 1) + (t + 2)) P_2) + 1,
\]
we can apply the inductive hypothesis, which states that the graph \( K_n \) contains a rainbow subgraph \( H \) such that \( H \cong (k - 1) P_4 \cup (t + 2) P_2 \).
Assume that \( H_1 \) and \( H_2 \) are two distinct subgraphs of \( H \), where \( H_1 = (k - 1) P_4 \) and \( H_2 = (t + 2) P_2 \). Furthermore, let \( H_1 = \{ P_4^i : i \in [k - 1] \} \) and \( H_2 = \{ P_2^j : j \in [t + 2] \} \). The edges of \( P_4^i \) can be defined as \( E(P_4^i) = \{ x_1^i x_2^i, x_2^i x_3^i, x_3^i x_4^i \} \) for \( i \in [k - 1] \), while the edges of \( P_2^j \) are given by \( E(P_2^j) = \{ y_1^j y_2^j \} \) for \( j \in [t + 2] \). 
Let \( G \) be a rainbow spanning subgraph of size \( |c(K_n)| \) derived from \( K_n \) that contains all edges of \( H \) (i.e., \( E(H) \subseteq E(G) \)). 
For the sake of contradiction, assume that \( K_n \) does not have a rainbow \( kP_4 \cup tP_2 \). From this assumption, we can make the following claims. To begin, we can derive Equation \eqref{eq-p2-p2} for this case using a method similar to that of Claim \ref{cl1-p2-p2}. For the sake of brevity, we will refer to this equation instead of restating it. Additionally, by applying a method analogous to Claim \ref{cl2-2n-4t-8}, we can establish the following claim.
\begin{claim}\label{cl3}
	It holds that
\begin{equation}\label{eqcl3}
|E(G[V'])|+|[V(H_2),V']_G|\leq2n - 4t - 8k.
\end{equation}
\end{claim}

Before starting the next claim, we need some definitions and notations. For \(h\) ranging from \(1\) to \(k-1\), we define the subset \(\mathcal{H}_h\) of \(H_1\) as follows:\\
1. The set \(\mathcal{H}_1\) is the largest subset of \(H_1\) such that for any two elements \(P^i_4\) and \(P^j_4\) in this set, the condition \(|E(G[V(P^i_4)\cup V(P^j_4)])| \geq 26\) holds.\\
2. Set \(F_1 = H_1\). For \(h\) ranging from \(2\) to \(k-1\), define \(F_h = F_{h-1} - \mathcal{H}_{h-1}\) and consider \(\mathcal{H}_h\) as the largest subset of \(F_h\) such that for any two elements \(P^i_4\) and \(P^j_4\) in this set, the condition \(|E(G[V(P^i_4)\cup V(P^j_4)])| \geq 26\) holds.

Additionally, let \(\iota\) be the largest integer in \([k-1]\) such that \(\mathcal{H}_\iota \neq \emptyset\). For \(h \in [\iota]\), define \(\mathcal{V}_h = \bigcup_{P \in \mathcal{H}_h} V(P)\).

\begin{claim}\label{cl5}
	For any \(h \in [\iota]\), it holds that  \(|[\mathcal{V}_h, V']_G| \leq3|V'|\).
\end{claim}
On the contrary, suppose \(|[V', \mathcal{V}_h]_G| \geq 3|V'| + 1\). In this case, since $n\geq8k+2t-4$, it follows that $|V'|\geq|\mathcal{V}_h|$. Thus, there exist four distinct vertices \(v_1, v_2, v_3, v_4 \in V'\) and four distinct vertices \(u_1, u_2, u_3, u_4 \in \mathcal{V}_h\) such that \(v_iu_i \in E(G)\) for \(i \in [4]\). Assume that \(u_1, u_2, u_3,\) and \(u_4\) belong to the \(r\) elements \(P_4^{i_1}, \ldots, P_4^{i_r}\) within \(\mathcal{H}_h\). Let \(V_1 = \{v_i : i \in [4]\}\) and \(U_1 = \{u_i : i \in [4]\}\). Based on our assumptions, for \(a \neq b \in [r]\), we have the following inequalities:

\begin{align}
  &|E(G[V(P_4^{i_a}) \cup V(P_4^{i_b})])| \geq 26,\label{cl5eq1}\\
  &|[V(P_4^{i_a}), V(P_4^{i_b})]_G| \geq 14.\label{cl5eq2}
\end{align}

To prove the claim, we will analyze four cases based on the value of \(r\). Notice that we will only investigate scenarios that cannot be obtained from each other by relabeling the vertices.  

\noindent
{\bf Case 1.} \(r = 1\). 
In this case, since \(G[U_1]\) contains a rainbow subgraph that is isomorphic to \(P_4\), it follows that \(G[U_1 \cup V_1]\) contains a rainbow subgraph isomorphic to \(2P_4\). However, this is sufficient to construct a rainbow subgraph of \(K_n\) that is isomorphic to \(kP_4 \cup tP_2\), leading to a contradiction.  \\
{\bf Case 2.} Let \( r = 2 \). In this case, we need to consider the following scenarios:\\
\noindent
{\bf Scenarios 2.1.} Let \( s \in [3] \) and \( z_1 \in ([4] - \{s, s+1\}) \) such that \( \{x_{s}^{i_1}, x_{s+1}^{i_1}, x_{z_1}^{i_1}\} \subseteq U_1 \). Assume that \( f \) is the element of \( [4] \) such that \( \{x^{i_2}_f\} = U_1 - \{x_{s}^{i_1}, x_{s+1}^{i_1}, x_{z_1}^{i_1}\} \). 
In this scenario, by employing \eqref{cl5eq2}, we have \( |[\{x_{z_1}^{i_1}, x_{z_2}^{i_1}\}, V(P^{i_2}_4)]_G| \geq 6 \), where \( z_2 \in ([4] - \{s, s+1, z_1\}) \).

 If \( |[\{z_1, z_2\}, V(P^{i_2}_4)]_G| = 6 \), then by using \eqref{cl5eq1}, we conclude that \( G[V(P^{i_2}_4)] \cong K_4 \).
Since \( |[\{z_1, z_2\}, V(P^{i_2}_4)]_G| = 6 \), it follows that there exist \( f_1 \neq f_2 \in ([4] - \{f\}) \) such that \( x_{z_1}^{i_1}x_{f_1}^{i_2}, x_{z_2}^{i_1}x_{f_2}^{i_2} \in E(G) \). Consequently, the three paths $v_{s}x_{s}^{i_1}x_{s+1}^{i_1}v_{s+1}$, $x_{z_2}^{i_1}x_{f_2}^{i_2}x_{f}^{i_2}v_{f}$, and $v_{z_1}x_{z_1}^{i_1}x_{f_1}^{i_2}x_{f_3}^{i_2}$, where \( f_3 \in ([4] - \{f, f_1, f_2\}) \), create a rainbow subgraph of \( G[V(P^{i_1}_4) \cup V(P^{i_2}_4) \cup V_1] \) that is isomorphic to \( 3P_4 \), which leads to a contradiction.

If \( |[\{z_1, z_2\}, V(P^{i_2}_4)]_G| = 7 \), then one of \( x_{z_1}^{i_1} \) and \( x_{z_2}^{i_1} \) has degree 4 in \( G[\{z_1, z_2\} \cup V(P^{i_2}_4)] \), while the other has degree greater than or equal to 3. Additionally, using \eqref{cl5eq1}, we find that \( G[V(P^{i_2}_4)] \) can be obtained from \( K_4 \) by removing an edge. Therefore, any three-subset of \( V(P^{i_2}_4) \) induces a subgraph of \( G \) that is either \( P_3 \) or \( K_3 \).
For \( g_1 \neq g_2 \in [4] \), suppose \( x_{g_1}^{i_2} \) is the farthest pendant vertex in \( P^{i_2}_4 \) from \( x_{f}^{i_2} \), and let \( x_{g_2}^{i_2} \) be the neighbor of \( x_{g_1}^{i_2} \) in \( P^{i_2}_4 \). It is evident that \( g_2 \neq f \). Based on earlier statements, for $g_3\in([4]-\{f,g_1,g_2\})$, the graph \( G[\{x_{g_1}^{i_2}, x_{g_2}^{i_2}, x_{g_3}^{i_2}, x_{z_2}^{i_1}\}] \) contains a rainbow subgraph \( L \) that is isomorphic to \( P_4 \). 
Therefore, if \( x_{z_1}^{i_1}x_{f}^{i_2} \in E(G) \), then the three paths $v_{s}x_{s}^{i_1}x_{s+1}^{i_1}v_{s+1}$, $v_{z_1}x_{z_1}^{i_1}x_{f}^{i_2}v_{f}$, and $L$ create a rainbow subgraph of \( G[V(P^{i_1}_4) \cup V(P^{i_2}_4) \cup V_1] \) that is isomorphic to \( 3P_4 \), resulting in a contradiction. If \( x_{z_1}^{i_1}x_{f}^{i_2} \notin E(G) \), then the three paths $v_{s}x_{s}^{i_1}x_{s+1}^{i_1}v_{s+1}$, $v_{z_1}x_{z_1}^{i_1}x_{g_1}^{i_2}x_{g_2}^{i_2}$, and $x_{z_2}^{i_1}x_{g_3}^{i_2}x_{f}^{i_2}v_{f}$ create a rainbow subgraph of \( G[V(P^{i_1}_4) \cup V(P^{i_2}_4) \cup V_1] \) that is also isomorphic to \( 3P_4 \), leading again to a contradiction. 

Otherwise, we have \( |[\{z_1, z_2\}, V(P^{i_2}_4)]_G| = 8 \). In this scenario, the three paths $v_{s}x_{s}^{i_1}x_{s+1}^{i_1}v_{s+1}$, $v_{z_1}x_{z_1}^{i_1}x_{g_1}^{i_2}x_{g_2}^{i_2}$, and $x_{z_2}^{i_1}x_{g_3}^{i_2}x_{f}^{i_2}v_{f}$ again create a rainbow subgraph of \( G[V(P^{i_1}_4) \cup V(P^{i_2}_4) \cup V_1] \) that is isomorphic to \( 3P_4 \), resulting once more in a contradiction.\\
\noindent
{\bf Scenarios 2.2.} Let \( s \), \( f_1 \), and \( f_2 \) be three elements in the set \([4]\) such that \( U_1 = \{ x_s^{i_1}, x_{s+1}^{i_1}, x_{f_1}^{i_2}, x_{f_2}^{i_2} \} \). In this scenario, by applying equation (2), for distinct \( z_1 \) and \( z_2 \) from the set \([4] - \{s, s+1\}\), we observe that \( |[\{x_{z_1}^{i_1}, x_{z_2}^{i_1}\}, V(P^{i_2}_4)]_G| \geq 6 \).

If \( |[\{x_{z_1}^{i_1}, x_{z_2}^{i_1}\}, V(P^{i_2}_4)]_G| = 6 \), then there exists \( g \in \{f_1, f_2\} \) such that \( x_s^{i_1} x_g^{i_2} \in E(G) \). Additionally, by using equation (1), we find that both \( G[V(P^{i_1}_4)] \) and \( G[V(P^{i_2}_4)] \) are isomorphic to \( K_4 \).
Suppose \( x_{g_1}^{i_2} \) is a neighbor of \( x_g^{i_2} \) in \( P^{i_2}_4 \). Define \( V(P^{i_2}_4) - \{x_g^{i_2}, x_{g_1}^{i_2}\}=\{x_{g_2}^{i_2}, x_{g_3}^{i_2}\} \). Under these conditions, we can see that there is \( h_1 \in \{g_2, g_3\} \) such that \( x_{z_1}^{i_1} x_{h_1}^{i_2} \in E(G) \). Let \( h_2 \in (\{g_2, g_3\} - \{h_1\}) \).
In this scenario, we can redefine the path \( P^{i_1}_4 \) to be \( x_{g_1}^{i_2} x_g^{i_2} x_s^{i_1} x_{s+1}^{i_1} \), and redefine the path \( P^{i_2}_4 \) as \( x_{z_2}^{i_1} x_{z_1}^{i_1} x_{h_1}^{i_2} x_{h_2}^{i_2} \). This allows us to apply either Case 1 or Scenario 2.1, leading us to reach a contradiction. 

If \( |[\{x_{z_1}^{i_1}, x_{z_2}^{i_1}\}, V(P^{i_2}_4)]_G| = 7 \), then by employing \eqref{cl5eq1}, either \( G[V(P^{i_1}_4)] \) or \( G[V(P^{i_2}_4)] \) is isomorphic to \( K_4 \). Therefore, the graph \( G[\{x_{z_1}^{i_1}, x_{z_2}^{i_1}, x_{g_2}^{i_2}, x_{g_2}^{i_2}\}] \) contains a rainbow subgraph \( L_1 \) that is isomorphic to \( P_4 \). Consequently, we can redefine the paths \( P^{i_1}_4 \) and \( P^{i_2}_4 \) as follows: \( P^{i_1}_4 \) becomes \( x_{g_1}^{i_2} x_g^{i_2} x_s^{i_1} x_{s+1}^{i_1} \) and \( P^{i_2}_4 \) becomes \( L_1 \). We can then apply either Case 1 or Scenario 2.1, leading us to a contradiction.

 If \( |[\{x_{z_1}^{i_1}, x_{z_2}^{i_1}\}, V(P^{i_2}_4)]_G| = 8 \), we can redefine the paths \( P^{i_1}_4 \) and \( P^{i_2}_4 \) as follows: \( P^{i_1}_4 \) will become \( x_{g_1}^{i_2} x_g^{i_2} x_s^{i_1} x_{s+1}^{i_1} \), and \( P^{i_2}_4 \) will become \( x_{z_1}^{i_1} x_{g_2}^{i_2} x_{z_2}^{i_1} x_{g_3}^{i_2} \). We can then apply either Case 1 or Scenario 2.1 to reach a contradiction.\\
 \noindent
{\bf Scenarios 2.3.} Let \(s^{i_1}_1\), \(s^{i_1}_2\), \(s^{i_2}_1\), and \(s^{i_2}_2\) be four elements in the set \([4]\) such that \(s^{i_1}_1 < s^{i_1}_2\), \(s^{i_2}_1 < s^{i_2}_2\), and 
$U_1 = \{ x^{i_1}_{s^{i_1}_1}, x^{i_1}_{s^{i_1}_2}, x^{i_2}_{s^{i_2}_1}, x^{i_2}_{s^{i_2}_2} \}$.
Additionally, let \([4] - \{s^{i_1}_1, s^{i_1}_2\} = \{f_1, f_2\}\) and \([4] - \{s^{i_2}_1, s^{i_2}_2\} = \{g_1, g_2\}\). Using the equations \eqref{cl5eq1} and \eqref{cl5eq2}, we have 
$|E(G[V(P_4^{i_1}) \cup V(P_4^{i_2})])| \geq 26$ and $|[V(P_4^{i_1}), V(P_4^{i_2})]_G| \geq 14$.

If we find that $|[V(P_4^{i_1}), V(P_4^{i_2})]_G| \leq 15$, then there exists \(l_1 \in \{i_1, i_2\}\) such that \(G[V(P_4^{l_1})] \cong K_4\). Suppose the remaining index is \(\{i_1, i_2\} - \{l_1\} = \{l_2\}\). Under these conditions, we can redefine \(P^{i_1}_4\) to become a rainbow subgraph of \(G[V(P^{l_1}_4)]\) that is isomorphic to \(P_4\) and contains the edge \(x^{l_1}_{s^{l_1}_1}x^{l_1}_{s^{l_1}_2}\). We consider the path \(P^{i_2}_4\) to be modified to become the path \(P^{l_2}_4\). We can then return to Scenario 2.2 to achieve our goal.

Alternatively, if $|[V(P_4^{i_1}), V(P_4^{i_2})]_G| = 16$, we will modify the path \(P^{i_1}_4\) to become \(x^{i_1}_{s^{i_1}_1}x^{i_2}_{s^{i_2}_1}x^{i_1}_{s^{i_1}_2}x^{i_2}_{s^{i_2}_2}\), and the path \(P^{i_2}_4\) will become \(x^{i_1}_{f_1}x^{i_2}_{g_1}x^{i_1}_{f_1}x^{i_2}_{g_2}\). We can then revert to Case 1 to achieve our desired outcome.\\

Before we begin proving other cases, we need to establish the following claim.

\begin{claim}\label{cl6}
Let \( l_1 \) and \( l_2 \) be two distinct elements in the set \( \{i_j:j\in[r]\} \). If \( V(P^{l_1}_4) \cap U_1 = \{u_h\} \) for some \( h \in [4] \), then there exists a vertex \( h_1 \in (V(P^{l_1}_4) - \{u_h\}) \) such that \( |[V(P^{l_2}_4), \{h_1\}]_G| = 4 \), and \( G[(V(P^{l_1}_4) - \{h_1\}) \cup \{v_h\}] \) contains a rainbow subgraph that is isomorphic to \( P_4 \).
\end{claim}

To prove this claim by contradiction, we assume it is false. This implies that for every vertex \( h_1 \in (V(P^{l_1}_4) - \{u_h\}) \), either \( |[V(P^{l_2}_4), \{h_1\}]_G| \leq 3 \) or \( G[(V(P^{l_1}_4) - \{h_1\}) \cup \{v_h\}] \) does not contain a rainbow subgraph that is isomorphic to \( P_4 \).
Let \( S \) denote the largest subset of \( V(P^{l_1}_4) - \{u_h\} \) such that for any \( s \in S \), we have \( |[V(P^{l_2}_4), \{s\}]_G| \leq 3 \). Additionally, define \( S_1 = V(P^{l_1}_4) - (\{u_h\} \cup S) \).

If \( |S| \geq 3 \), we can conclude that \( |E(G[V(P_4^{l_1}) \cup V(P_4^{l_2})])| \leq 25 \), which contradicts equation \eqref{cl5eq1}. 

Now, consider the case where \( |S| \leq 2 \). Since for any element \( s_1 \in S_1 \), the graph \( G[(V(P^{l_1}_4) - \{s_1\}) \cup \{v_h\}] \) does not contain a rainbow subgraph isomorphic to \( P_4 \), it follows that for each pair of vertices \( s_1, s_1' \in S_1 \), either \( s_1 s_1' \notin E(G) \) or \( s_1 u_h, s_1' u_h \notin E(G) \). Therefore, we have:
\[
|E(G[S_1 \cup \{u_h\}])| \leq \binom{|S_1| + 1}{2} - \binom{|S_1|}{2}.
\]
Consequently, since \( |S| = 4 - |S_1| \), we can compute:

\begin{align*}
|E(G[V(P_4^{l_1}) \cup V(P_4^{l_2})])| \leq & 28 - \left(4 - |S_1| + \binom{|S_1|}{2}\right) \\
= & 28 - \left(4 + \frac{|S_1|(|S_1| - 3)}{2}\right) \leq 25, \quad \text{since } |S_1| \in [4].
\end{align*}

This again leads to a contradiction with equation \eqref{cl5eq1}. Therefore, we conclude that our claim is true. \qed\\

\noindent
{\bf Case 3.}  Let \( r = 3 \). In this case, there are three distinct elements \( l_1, l_2, l_3 \in \{i_1,i_2,i_3\} \) such that \( |V(P^{l_1}_4) \cap U_1| = 2 \) and \( |V(P^{l_2}_4) \cap U_1| = |V(P^{l_3}_4) \cap U_1| = 1 \). 
Suppose \( s \) and \( f \) are two elements from the set \( [4] \) such that \( V(P^{l_2}_4) \cap U_1 = \{u_s\} \) and \( V(P^{l_3}_4) \cap U_1 = \{u_f\} \). 
Using Claim \ref{cl6}, we know there exist elements \( s_1 \in (V(P^{l_2}_4) - \{u_s\}) \) and \( f_1 \in (V(P^{l_3}_4) - \{u_f\}) \) such that \( |[\{s_1\}, V(P_4^{l_1})]_G| = |[\{f_1\}, V(P_4^{l_1})]_G| = 4 \). Additionally, each of the graphs \( G[(V(P_4^{l_2}) - \{s_1\}) \cup \{v_s\}] \) and \( G[(V(P_4^{l_3}) - \{f_1\}) \cup \{v_f\}] \) contains a rainbow subgraph that is isomorphic to \( P_4 \). 
For \( g \neq g' \in ([4] - \{s, f\}) \), let \( V(P^{l_1}_4) \cap U_1 = \{u_g, u_{g'}\} \) and \( V(P^{l_1}_4) - \{u_g, u_{g'}\} = \{a_1, a_2\} \). Under these conditions, we can see that \( a_1 s_1 u_g v_g \) and \( a_2 f_1 u_{g'} v_{g'} \) are two distinct rainbow subgraphs of \( G[V(P_4^{l_1}) \cup \{s_1, f_1\}] \), and each of them is isomorphic to \( P_4 \). This leads us to a contradiction.

\noindent
{\bf Case 4.}  Let \( r = 4 \). In this scenario, for each \( l \in [4] \), we have \( |V(P^{i_l}_4) \cap U_1| = 1 \). For \( l \in [3] \), let \( V(P^{i_l}_4) \cap U_1 = \{u_{l_1}\} \). By employing Claim \ref{cl6}, for each \( l \in [3] \), there exists an \( f_l \in (V(P^{i_l}_4) - \{u_{l_1}\}) \) such that \( |[\{f_l\}, V(P_4^{i_4})]_G| = 4 \). Furthermore, for each \( l \in [3] \), the graph \( G[(V(P_4^{i_l}) - \{f_l\}) \cup \{v_{l_1}\}] \) contains a rainbow subgraph that is isomorphic to \( P_4 \). 
Assume \( V(P_4^{i_4}) \cap U_1 = \{u_g\} \), where \( g \in ([4] - \{l_1 : l \in [3]\}) \), and assume \( V(P_4^{i_4}) - \{u_g\} = \{a_1, a_2, a_3\} \). Under these conditions, we can observe that \( v_g u_g f_1 a_1 \) and \( a_2 f_2 a_3 f_3 \) are two distinct rainbow subgraphs of \( G[V(P_4^{i_4}) \cup \{f_i : i \in [3]\}] \), and each is isomorphic to \( P_4 \). This leads us to a contradiction.

From the four cases discussed, we conclude that Claim \ref{cl5} is true. \qed\\

Let \(\Gamma = |E(G[V(H_1)])| + |[V(H_1), V']_G|\). Let \(i \in ([\iota] - \{1\})\) and \(j \in [i-1]\). By considering the maximality of \(|\mathcal{H}_j|\), we can observe that for each \(F \in \mathcal{H}_j\), there exists at least one \(F' \in \mathcal{H}_i\) such that \(|E(G[V(F) \cup V(F')])| \leq 25\). 
Furthermore, any edges in the set \(E(G[V(H_1)])\) can belong to at most \(\iota - 1\) graphs of the form \(\mathcal{H}_j \cup \mathcal{H}_i\), where \(i \neq j \in [\iota]\). 
Therefore, by applying Claim \ref{cl5}, we can conclude that:
\begin{align*}
\Gamma &\leq \frac{1}{2}(4k-4)(4k-5) + \sum_{i=1}^{\iota}|[\mathcal{V}_i, V']_G|  -\frac{1}{\iota-1} \sum_{i=1}^{\iota - 1} \sum_{j=i + 1}^{\iota} 3|\mathcal{H}_j| \\
&\leq  \frac{1}{2}(4k-4)(4k-5) + 3|V'|\iota - \frac{3}{2}\iota.
\end{align*}
Let us define the function \(f(\iota) = 3|V'|\iota - \frac{3}{2}\iota\). Since \(\iota \leq k - 1 \leq |V'|\), it follows that the function \(f\) is increasing on the interval \([1, k - 1]\). Therefore, from this and the previous inequality, we can derive the following inequality:

\begin{equation}\label{eqcl6}
\Gamma \leq \frac{1}{2}(4k - 4)(4k - 5) + 3(k - 1)(n - 4k - 2t)-\frac{3}{2}(k-1).
\end{equation}

\begin{claim}\label{cl-h2-h1}
Let $\Upsilon=|[V(H_1), V(H_2)]_G|$. It holds that 
\begin{equation}\label{eqcl7}
\Upsilon\leq2(k-1)(t+3).
\end{equation}
\end{claim}

To prove this claim, consider \( i \in [k-1] \) and distinct indices \( z_1, z_2\in [t+2] \). If there exist values \( h_1, h_2 \in [2] \), \( r_1 \in [2] \), and \( r_2 \in ([4] - [2]) \), such that \( y^{z_1}_{h_1} x^i_{r_1} \in E(G) \) and \( y^{z_2}_{h_2} x^i_{r_2} \in E(G) \), then each of the subgraphs \( G[V(P_2^{z_1}) \cup \{x^i_1, x^i_2\}] \) and \( G[V(P_2^{z_2}) \cup \{x^i_3, x^i_4\}] \) contains a rainbow subgraph that is isomorphic to \( P_4 \). Consequently, the graph \( G[V(P_4^i) \cup  V(P_2^{z_1})\cup V(P_2^{z_2})] \) will include a rainbow subgraph that is isomorphic to \( 2P_4 \), which contradicts our assumptions.
Therefore, for every pair \( z_1 \neq z_2 \in [t+2] \), and for any \( h_1, h_2 \in [2] \), \( r_1 \in [2] \), and \( r_2 \in ([4] - [2]) \), at least one of the following conditions must hold: either \( y^{z_1}_{h_1} x^i_{r_1} \not\in E(G) \) or \( y^{z_2}_{h_2} x^i_{r_2} \not\in E(G) \). 
Moreover, by similar reasoning, we can conclude that 
\[
|[\{x^i_1, x^i_2\}, \cup_{j=1}^{t+2} V(P_2^{j})]_G| \leq 2t+6 \quad \text{and} \quad |[\{x^i_3, x^i_4\}, \cup_{j=1}^{t+2} V(P_2^{j})]_G| \leq 2t+6.
\]
As a result, the size of the set \( [V(P_4^i), \cup_{j=1}^{t+2} V(P_2^{j})]_G \) can be at most $2t+6$. Thus, we can conclude that 

\[
\Upsilon \leq (k-1)(2t+6)= 2(k-1)(t+3).
\]
This concludes the proof of the claim. \qed\\

Now, by utilizing Equations \eqref{eq-p2-p2}, \eqref{eqcl3}, \eqref{eqcl6}, and \eqref{eqcl7}, we conclude that 

\begin{equation}\label{eq8}
|c(K_n)| = |E(G)| \leq \frac{1}{2} \left( 6kn+ 2t + 15 -8k^2  - 8kt - 19k - 2n  \right).
\end{equation}

Using our assumption, we also have that \( |c(K_n)| = AR(n, (2k+t)P_2) + 1 \). Since \( n \geq 8k + 2t - 4 \) and \( k \geq 2 \), applying Equation \eqref{eqth0} gives us the following inequality:

\begin{equation}\label{eq9}
|c(K_n)|= 
\begin{cases} 
(2k + t - 2)(4k + 2t - 3) + 2 & \text{when } n \leq \frac{5(2k+t)-7}{2}, \\ 
(2k + t - 2)\left(n - \frac{2k+t-1}{2}\right) + 2 & \text{when } n \geq \frac{5(2k+t)-7}{2}. 
\end{cases}
\end{equation}

We define the following:
\begin{align*}
\Theta & = \frac{1}{2} \left( 6kn+ 2t + 15 -8k^2  - 8kt - 19k - 2n  \right), \\
\Phi & = (2k + t - 2)(4k + 2t - 3) + 2, \\
\Psi & = (2k + t - 2)\left(n - \frac{2k+t-1}{2}\right) + 2.
\end{align*}

If \( n \leq \frac{5(2k+t)-7}{2} \), then we conclude that:

\begin{align*}
\Phi - \Theta &\geq\frac{1}{2}(24k^2+ 24kt + 4t^2+2n+1  - 6kn - 9k  - 16t)\\
&\geq \frac{1}{2}( 9kt + 4t^2 + 22k-6k^2  - 11t - 6),{\rm as }~k\geq2~{\rm and~} n \leq \frac{5(2k+t)-7}{2}.
\end{align*} 
So, by employing our assumption $t\geq k+1$, we can observe that
\[\Phi - \Theta \geq\frac{1}{2}(7k^2 + 28k - 13)>0,\]
which creates a contradiction between Equations \eqref{eq8} and \eqref{eq9}.  

If \( n \geq \frac{5(2k+t)-7}{2} \), then we can check that:
\begin{align*}
 \Psi-\Theta&\geq\frac{1}{2}(4k^2+ 4kt + 2nt+ 25k  + t- 2kn  - t^2  - 2n - 13)\\
 &\geq \frac{1}{2}(9kt + 4t^2 + 22k-6k^2  - 11t - 6),{\rm as }~t\geq k+1~{\rm and~} n \geq \frac{5(2k+t)-7}{2}.
\end{align*}
So, since $t\geq k+1$, it follows that:
\[
\Psi - \Theta \geq \frac{1}{2}(7k^2 + 28k - 13) > 0,
\]
which again creates a contradiction between Equations \eqref{eq8} and \eqref{eq9}.

Thus, from the above arguments, we conclude that Equation \eqref{eqth4} holds true, thereby completing the proof of Theorem \ref{th1}. \qed\\

The proof of Theorem~\ref{th2} closely resembles the proof of Theorem \ref{th1}. We simply need to consider \( V' = \emptyset \) and then follow the same steps as in the proof of Theorem \ref{th1} to achieve our objective. Therefore, we will omit the details here for brevity.

\section{Proof of Theorem~\ref{th3}}
For two positive integers \( t \) and \( n \), assume that \( t \geq 1 \) and \( n \geq 2t + 4 \). When \( t = 1 \), we can use Equations \eqref{eqth1} and \eqref{eqth2} to demonstrate that for \( n \geq 6 \), \( AR(n, P_4 \cup P_2) = AR(n, 3P_2) = n \). This confirms the theorem when \( t = 1 \). 
Additionally, by employing Equation \eqref{eqth3}, we have 
\begin{equation*}
AR(n, P_4 \cup tP_2) \geq AR(n, (t+2)P_2).
\end{equation*}

To establish our result, we need to verify the following inequality for \( t \geq 2 \):
\begin{equation}\label{lleq1}
AR(n, P_4 \cup tP_2) \leq AR(n, (t+2)P_2).
\end{equation}

To prove this, we must show that if \( c \) is an arbitrary edge-coloring of \( K_n \) such that \( |c(K_n)| = AR(n, (t+2)P_2) + 1 \), then \( K_n \) contains a rainbow subgraph isomorphic to \( P_4 \cup tP_2 \). Since \( |c(K_n)| = AR(n, (t+2)P_2) + 1 \), it follows that there exists a rainbow subgraph \( H \) of \( K_n \) such that \( H \cong (t+2)P_2 \). We will divide the remainder of the proof into two scenarios:

\noindent
{\bf Scenario 2.1:} \( n = 2t + 4 \). In this scenario, by applying Equation \eqref{eqth1}, we can observe that 
\[
|c(K_n) - c(H)|  = \begin{cases} 
	\frac{1}{2}t(3t+5) & {\rm~when~ } 1\leq t\leq 4, \\ 
	2t^2+1 & {\rm~when~ } t\geq 5. 
\end{cases}
\]
Since $t\geq2$, it follows that $c(K_n) - c(H)\neq\emptyset$. Therefore, using \( H \) and an edge \( e \) from \( K_n \) such that \( c(e) \in (c(K_n) - c(H)) \), we can construct a rainbow subgraph of $K_n$ that is isomorphic to \( P_4 \cup tP_2 \).

\noindent
{\bf Scenario 2.1:} Let \( n \geq 2t + 5 \). In this case, using Equation \eqref{eqth0}, we find that:
\[
|c(K_n)| = \begin{cases} 
	t(2t + 1) + 2 & \text{when } n \leq \frac{5t + 3}{2}, \\ 
	t(n - \frac{t + 1}{2}) + 2 & \text{when } n \geq \frac{5t + 3}{2}. 
\end{cases}
\]
Furthermore, if we assume that \( K_n \) does not contain a rainbow \( P_4 \cup tP_2 \), then by applying a method similar to that used in the case \( k = 1 \) in the proof of Theorem \ref{th1}, we can conclude that:
\[
|c(K_n)| \leq 2n - 3t - 6.
\]

Next, we will analyze the two cases separately: \( n \leq \frac{5t + 3}{2} \) and \( n \geq \frac{5t + 3}{2} \). Since \( t \geq 2 \), this leads us to the following conclusions:
\[
|c(K_n)| \leq \begin{cases} 
	2n - 3t - 6 \leq 2t - 3 < t(2t + 1) + 2 & \text{when } n \leq \frac{5t + 3}{2}, \\ 
	2n - 3t - 6 < t\left(n - \frac{t + 1}{2}\right) + 2 & \text{when } n \geq \frac{5t + 3}{2}. 
\end{cases}
\]
This leads to a contradiction with our assumption about \(|c(K_n)|\).

Based on these two scenarios, we conclude that if \( t \geq 2 \), then Equation \eqref{lleq1} holds true, which completes the proof of Theorem \ref{th3}. \qed

\section{Proof of Theorem~\ref{kp4}}
Let \( k \) and \( n \) be two positive integers such that \( k \geq 2 \) and \( n = 4k \). To prove our theorem, we need to confirm the following two inequalities:
\begin{align}
AR(n, kP_4) & \geq (2k-1)(4k-3) + 1, \label{eqkp4-l} \\
AR(n, kP_4) & \leq (2k-1)(4k-3) + 1. \label{eqkp4-u}
\end{align}

First, we will prove the lower bound. To do this, consider a subgraph \( K_{4k-2} \) of \( K_n \) and color it using a rainbow coloring. The remaining edges will be colored with one additional color. This construction results in an edge-coloring of \( K_n \) with \( (2k-1)(4k-3) + 1 \) colors, ensuring that there is no rainbow \( kP_4 \). This satisfies Equation \eqref{eqkp4-l}.

Next, we will prove the upper bound. We need to show that if \( c \) is an arbitrary edge-coloring of \( K_n \) such that \( |c(K_n)| = (2k-1)(4k-3) + 2 \), then \( K_n \) must contain a rainbow subgraph isomorphic to \( kP_4 \). For the sake of contradiction, assume that \( K_n \) does not have a rainbow subgraph isomorphic to \( kP_4 \). Since \( k \geq 2 \), it follows that:
\[
|c(K_n)| = (2k-1)(4k-3) + 2 \geq 
\begin{cases} 
\frac{1}{2}(2k-2)(6k+1) + 2 & \text{when } 2 \leq k \leq 3, \\ 
(2k-2)(4k-3) + 3 & \text{when } k \geq 4. 
\end{cases}
\]
Therefore, by employing Theorem \ref{th2} and Equation \eqref{eqth1}, we can conclude that 
\[
|c(K_n)| \geq AR(n, (k-1)P_4 \cup 2P_2) + 1.
\]
This implies that \( K_n \) has a rainbow subgraph \( H \) that is isomorphic to \( (k-1)P_4 \cup 2P_2 \). Assume $H_1$ and $H_2$ are two distinct subgraphs of $H$ that $H_1=(k-1)P_4$ and $H_2=2P_2$. Also, Let \( G \) be a rainbow spanning subgraph of size \( |c(K_n)| \) from \( K_n \) such that \( E(H) \subseteq E(G) \).
It is straightforward to see that \( c(G[V(H_2)]) - c(H_2) = \emptyset \), as otherwise we could find a rainbow \( kP_4 \), which would lead to a contradiction. Therefore, we have \( |E(G[V(H_2)])| = 2 \). By employing a similar method as in Claim \ref{cl-h2-h1}, we can derive that 
\[
|[V(H_1), V(H_2)]_G| \leq 8(k-1), 
\]
noting that here \( t = 1 \). Additionally, it is evident that 
\[
|E(G[V(H_1)])| \leq \binom{4(k-1)}{2}. 
\]
Consequently, we derive:
\begin{align*}
|E(G)| &= |E(G[V(H_1)])| + |E(G[V(H_2)])| + |[V(H_2), V(H_1)]_G| \\ 
&\leq \binom{4(k-1)}{2} + 2 + 8(k-1) \\ 
&= (2k-1)(4k-3) + 1.
\end{align*}
This leads us to a contradiction with our assumption regarding the size of \( c(K_n) \). Thus, if \( |c(K_n)| = (2k-1)(4k-3) + 2 \), then \( K_n \) must contain a rainbow subgraph that is isomorphic to \( kP_4 \). Therefore, Equation \eqref{eqkp4-u} holds true. This completes the proof of Theorem \ref{kp4}. \qed
\section{Concluding remarks}
Let \(G = kP_4 \cup tP_2\) for three non-negative integers \(k\), \(t\), and \(n\), where \(n \geq 4k + 2t\). We define \(E'\) as a subset of the edge set \(E(G)\) such that each endpoint of these edges has a degree of two in \(G\). In this study, we explored the anti-Ramsey number in the context of edge deletion, demonstrating that both decreasing and unchanging behaviors are possible outcomes.
We prove that the behavior of the anti-Ramsey number remains consistent when the edges in \(E'\) are removed from \(G\), i.e., \(AR(n, G) = AR(n, G - E')\), if one of the conditions: (i) \(t \geq k + 1 \geq 2\) and \(n \geq 8k + 2t - 4\); (ii) \(k, t \geq 1\) and \(n = 4k + 2t\); and (iii) \(k = 1\), \(t \geq 1\), and \(n \geq 2t + 4\) holds. However, this consistency does not hold when \(k \geq 2\), \(t = 0\), and $n=4k$. As a result, we calculate \(AR(n,kP_4 \cup tP_2)\) for the cases: (i) \(t \geq k + 1 \geq 2\) and \(n \geq 8k + 2t - 4\); (ii) \(k, t \geq 1\) and \(n = 4k + 2t\); 
(iii) \(k = 1\), \(t \geq 0\), and \(n \geq 2t + 4\); and (iv) \(k \geq 1\), \(t = 0\), and \(n = 4k\). For future research, the following questions may be considered:
\begin{problem} 
Let \(G\) be a graph with an order of at least \(n\). 
\begin{enumerate}
    \item  Is there a subset \(S\) of \(E(G)\) such that \(|S| \geq 1\) and \(AR(n, G) = AR(n, G - S)\)?
    \item  What is the size of the largest subset \(S\) of \(E(G)\) such that: 
    \[AR(n, G) = AR(n, G - S)? \]
    \item For which positive integers \( k \), \( t \), \( a \), \( b \), and \(n\) with the property \( a = 2b \) and \(n\geq ka+tb\) does it hold that: 
\[ AR(n, kP_a \cup tP_b) = AR(n, (2k + t)P_b)? \]
\end{enumerate}
\end{problem}

\vskip3mm

\noindent \textbf{Acknowledgments}
\vspace{2mm}

The research of Ali Ghalavand and Xueliang Li was supported the by the Natural Science Foundation of China under Grant No.12131013.  Jin was supported by the Natural Science Foundation of China under Grant No.12571380.






\begin{thebibliography}{99}

\bibitem{R2-SP-1} S. Akbari, A. Alipour, Multicolored trees in complete graphs, J. Graph Theory, 54 (2007), 221-232.

\bibitem{p1-1} A. Bialostocki, S. Gilboa, Y. Roditty, Anti-Ramsey numbers of small graphs, Ars Combin. 123 (2015), 41-53.

\bibitem{R2-SP-2} A. Bialostocki, W. Voxman, On the anti-Ramsey numbers for spanning trees, Bull. Inst.Combin. Appl. 32 (2001),  23-26.

\bibitem{p1-2} H. Chen, X.L. Li, J.H. Tu, Complete solution for the rainbow number of matchings, Discrete Math. 309 (10) (2009), 3370-3380.

 \bibitem{Er1} P. Erd\H{o}s, M. Simonovits, V.T. S\'{o}s, Anti-Ramsey theorems, Colloq. Math. Soc. J\'{a}nos Bolyai, Vol. 10 North-Holland Publishing Co., Amsterdam-London, 1975, pp. 633-
643.

\bibitem{p1-4} C.Q. Fang, E. Gy\H{o}ri, M. Lu, J.M. Xiao, On the anti-Ramsey number of forests, Discrete Appl. Math. 291 (11) (2021), 129-142.

\bibitem{p1-5} S. Fujita, A. Kaneko, I. Schiermeyer, K. Suzuki, A rainbow $k$-matching in the complete graph with $r$ colors, Electron. J. Combin. 16 (2009), \#R51.

\bibitem{GL-1} A. Ghalavand, X. Li, On the Anti-Ramsey Number of Spanning Linear Forests with Paths of Lengths 2 and 3, arXiv: 2509.25949 [math.CO] (30 Sep 2025).

\bibitem{p1-6} S. Gilboa, Y. Roditty, Anti-Ramsey numbers of graphs with small connected components, Graphs Combin. 32 (2) (2016), 649-662.

\bibitem{R2-HP-1} R. Gu, J. Li, Y. Shi, Anti-Ramsey numbers of paths and cycles in hypergraphs, SIAM J. Discrete Math. 34 (1) (2020), 271-307.

\bibitem{R2-HP-2} M. Guo, H. Lu, X. Peng, Anti-Ramsey number of matchings in 3-uniform hypergraphs, SIAM J. Discrete Math. 37 (3) (2023), 1970-1987.

\bibitem{p1-7} R. Haas, M. Young, The anti-Ramsey number of perfect matching, Discrete Math. 312 (2012), 933-937.

\bibitem{p1-8} M.L. He, Z.M. Jin, Rainbow short linear forests in edge-colored complete graph, Discrete Appl. Math. 361 (2025), 523-536.

\bibitem{p1-9} S. Jahanbekam, D.B. West, Anti-Ramsey problems for t edge-disjoint rainbow spanning subgraphs: cycles, matchings, or trees, J. Graph Theory 82 (1) (2016), 75-89.

\bibitem{new-3} Q. Jie, Z.M. Jin, Anti-Ramsey number of union of 5-path and matching, Discuss.  Math. Graph Theory 45 (2025), 1185-1210.

\bibitem{p1-10} Q. Jie, M.L. He, Z.M. Jin, Rainbow forest consisting of short paths in $K_n$, Discrete Appl. Math. 376 (2025), 260-269.

\bibitem{new-2} Q. Jie, Z.M. Jin, Rainbow-free colorings for spanning linear forest consisting of short paths, 2025, {\it submitted}.

\bibitem{new-1} Z.M. Jin and J.Q. Gu, Rainbow disjoint union of clique and matching in edge-colored complete graph, Discuss. Math. Graph Theory 44 (2024), 953–970.

\bibitem{n-p4-1} Z. Jin, Q. Jie, Z. Cao, Rainbow disjoint union of $P_4$ and a matching in complete graphs, Appl. Math. Comput. 474 (2024), 128679.

\bibitem{R2-SP-3} L.Y. Lu, A. Meier, Z.Y. Wang, Anti-Ramsey number of edge-disjoint rainbow spanning trees in all graphs, SIAM J. Discrete Math. 37 (2) (2023), 1162-1172.

\bibitem{R2-SP-4} L.Y. Lu, Z.Y. Wang, Anti-Ramsey number of edge-disjoint rainbow spanning trees, SIAM J. Discrete Math. 34 (4) (2020), 2346--2362.

\bibitem{p3-6} J.J. Montellano-Ballesteros, V. Neumann-Lara, An anti-Ramsey theorem on cycles, Graphs Combin. 21 (3) (2005), 343-354.

 \bibitem{p2-9} I. Schiermeyer, Rainbow numbers for matchings and complete graphs, Discrete Math. 286 (2004), 157-162.

\bibitem{p2-10} M. Simonovits, V.T. Sós, On restricted coloring of $K_n$, Combinatorica 4 (1) (1984), 101–110.

\bibitem{R2-HP-3} Y. Tang, T. Li, G. Yan, Anti-Ramsey number of expansions of paths and cycles in uniform hypergraphs, J. Graph Theory 101 (4) (2022), 668-685. 

\bibitem{p1-14} T.Y. Xie, L.T. Yuan, On the anti-Ramsey numbers of linear forests, Discrete Math. 343(12)(2020), 112130.


\end{thebibliography}
\end{document}